\documentclass[12pt,reqno]{amsart}
\usepackage{amsmath,amssymb,amsfonts,amscd,latexsym,amsthm,mathrsfs}
\usepackage{enumerate}
\usepackage{graphicx}
\usepackage[usenames]{color}
\usepackage[unicode]{hyperref}

\usepackage{cite}

\let\lf\lfloor
\let\rf\rfloor
\let\<\langle
\let\>\rangle

\newcommand{\qbin}[2]{\genfrac{[}{]}{0pt}{}{#1}{#2}}
\renewcommand{\pmod}[1]{\;\allowbreak(\operatorname{mod}#1)}
\newcommand{\KR}{\operatorname{KR}}
\newcommand{\RK}{\operatorname{\rotatebox[origin=c]{180}{$\operatorname{RK}$}}}

\begin{document}

\title[Reflecting (on) the modulo 9 Kanade--Russell identities]{Reflecting (on) the modulo 9\\Kanade--Russell (conjectural) identities}

\author{Ali Uncu}
\address{[A.U.] University of Bath, Faculty of Science, Department of Computer Science, Bath, BA2\,7AY, UK}
\email{aku21@bath.ac.uk}

\address{[A.U.] Johann Radon Institute for Computational and Applied Mathematics, Austrian Academy of Science, Altenbergerstraße 69, A-4040 Linz, Austria}
\email{akuncu@ricam.oeaw.ac.at}

\author{Wadim Zudilin}
\address{[W.Z.] Department of Mathematics, IMAPP, Radboud University, PO Box 9010, 6500~GL Nijmegen, Netherlands}
\email{w.zudilin@math.ru.nl}

\dedicatory{To Doron Zeilberger, with Experimental Mathematics wishes,\\on his twentieth prime birthday}

\thanks{Research of the first author is supported partly by EPSRC grant number EP/T015713/1 and partly by FWF grant P-34501-N}



\date{21 February 2021}

\begin{abstract}
We examine complexity and versatility of five modulo 9 Kanade--Russell identities through their finite (\emph{aka} polynomial) versions and images under the $q\mapsto1/q$ reflection.
\end{abstract}

\maketitle

\section{Introduction}
\label{sec1}

Every second paper about partition identities features the Rogers--Ramanujan identities as toy examples.
Let us follow this tradition and consider one of those burdened with fame as the limiting case of the \emph{finite} version
\begin{equation}
\sum_{n\ge0}q^{n^2}\qbin{N}{n}
=\sum_{k\in\mathbb Z}(-1)^kq^{k(5k+1)/2}\qbin{2N}{N+2k}
\label{rama1fin}
\end{equation}
due to Bressoud~\cite{Br81}.
Here and in what follows we make use of standard $q$-hypergeometric notation:
$$
\qbin Nn=\qbin{N}{n}_q=\begin{cases}
\dfrac{(q;q)_N}{(q;q)_n(q;q)_{N-n}} &\text{for}\; n=0,1,\dots,N, \\[1.5mm]
0 &\text{otherwise},
\end{cases}
$$
denotes a $q$-binomial coefficient,
$$
(a;q)_n=\prod_{j=0}^{n-1}(1-aq^j)
$$
is a $q$-shifted factorial ($q$-Pochhammer symbol), also meaningful for $n=\infty$, and
$$
(a_1,a_2,\dots,a_k;q)_n=(a_1;q)_n(a_2;q)_n\dotsb(a_k;q)_n.
$$
When $|q|<1$, the limit as $N\to\infty$ in \eqref{rama1fin} translates the equality of two polynomials into
\begin{equation}
\sum_{n\ge0}\frac{q^{n^2}}{(q;q)_n}
=\frac1{(q;q)_\infty}\sum_{k\in\mathbb Z}(-1)^kq^{k(5k+1)/2}
=\frac1{(q,q^4;q^5)_\infty},
\label{rama1}
\end{equation}
where Jacobi's triple product identity is applied to the sums on the right-hand sides.
On the other hand, we can \emph{reflect} the identity \eqref{rama1fin} in a different one by applying the involution $q\mapsto1/q$.
Since
\[\qbin{n+m}{m}_{q^{-1}} = q^{-nm} \qbin{n+m}{m}_q, \]
we find out that
\begin{equation*}
q^{-N^2}\sum_{n\ge0}q^{nN}\qbin{N}{n}
=q^{-N^2}\sum_{k\in\mathbb Z}(-1)^kq^{k(3k-1)/2}\qbin{2N}{N+2k}
\end{equation*}
(in the sum on the left-hand side we change the summation index $n\mapsto N-n$).
Apparently, this identity is not as interesting as the original \eqref{rama1fin}: the polynomial obtained from the multiplication of either side by $q^{N^2}$ tends to (boring) 1 as $N\to\nobreak\infty$.

Let us take a look at a contemporary and more advanced (in particular, conjectural!) example that comes out from the inspiring work \cite{KR15} of Kanade and Russell.
Denote by $s_N(n)$ the number of partitions $\lambda=(\lambda_1,\lambda_2,\dots)$ of $n$ such that
\begin{enumerate}[(a)]
\item \label{cond-a}
the largest part of $\lambda$ is at most $N$;
\item \label{cond-b}
$\lambda$ has no parts of size 1;
\item \label{cond-c}
the difference of parts at distance 2 is at least 3;
\item \label{cond-d}
if consecutive parts differ by at most 1 then their sum is congruent to $2\pmod3$.
\end{enumerate}
Then the generating function $S_N(q)=\sum_{n\ge0}s_N(n)q^n$ is (clearly) a polynomial for each $N\ge0$, and the polynomials $S_N(q)$ satisfy the recursion
\begin{equation*}
S_N(q)=S_{N-1}(q)+q^N\times\begin{cases}
(1+q^{N-1})S_{N-3}(q)+q^{N-2}S_{N-4}(q) &\text{if}\; N\equiv0\pmod3, \\
S_{N-2}(q)+q^NS_{N-3}(q) &\text{if}\; N\equiv1\pmod3, \\
S_{N-2}(q) &\text{if}\; N\equiv2\pmod3, \\
\end{cases}
\end{equation*}
for $N=1,2,\dots$, with initial conditions $S_0(q)=1$ and $S_N(q)=0$ for $N<0$.
The generating function $S_\infty(q)=\lim_{N\to\infty}S_N(q)$ counts partitions with the condition on the largest part dropped, that is, satisfying hypotheses \eqref{cond-b}--\eqref{cond-d} above, and the conjecture from~\cite{KR15} predicts that
\begin{equation}
S_\infty(q)\overset?=\frac{1}{(q^2;q^3)_\infty(q^3;q^9)_\infty}
=\frac{1}{(q^2,q^3,q^5,q^8;q^9)_\infty}.
\label{S-lim}
\end{equation}
As we will see later, the series $S_\infty(q)$ can be written as a double-sum Rogers--Ramanujan type identity\,---\,this was found by Kur\c sung\"oz \cite{Ku19}.
The conjecture in~\eqref{S-lim} has a surprisingly different (and supposedly more difficult) counterpart, which was observed by Warnaar~\cite{Wa16}:
\begin{align}
\lim_{M\to\infty} q^{M(3M+2)} S_{3M}(q^{-1})
&\overset?=\frac{1}{(q^2;q^3)_\infty (q^3,q^9,q^{12},q^{21},q^{30},q^{36},q^{39};q^{45})_\infty},
\nonumber\\
\lim_{M\to\infty} q^{M(3M+5)} S_{3M+1}(q^{-1})
&\overset?=\frac{1}{(q^2;q^3)_\infty (q^3,q^{12},q^{18},q^{21},q^{27},q^{30},q^{39};q^{45})_\infty},
\nonumber\\
\lim_{M\to\infty} q^{(M+1)(3M+2)} S_{3M+2}(q^{-1})
&\overset?=\lim_{M\to\infty} q^{M(3M+2)} S_{3M}(q^{-1})
\nonumber\\ &\quad
+q^2\lim_{M\to\infty} q^{M(3M+5)} S_{3M+1}(q^{-1}).
\label{conj:OW3}
\end{align}
One goal of this work is to analyze this $q\mapsto1/q$ phenomenon from a more general perspective, in particular, to provide other modulo~45 product sides of all other (five in total) Kanade--Russell modulo~9 partition counting functions.
On this way, we give explicit finite versions of the functions (for example, of $S_N(q)$ above) and explicit (Rogers--Ramanujan type) sums for the reflected identities.
One pleasing outcome of this routine is a \emph{proof} of~\eqref{conj:OW3} as well as other similar cases.

\medskip
Doron Zeilberger's combined interests in combinatorics of generalized Rogers--Ramanujan identities \cite{BZ89}, algorithmic aspects of $q$-identities and Experimental Mathematics have served a fruitful motivation for many, junior and senior, to produce beautiful research pieces.
One notable recent illustration is the discovery of novel identities for the classical partitions by Doron's mathematical descendant, Matthew Russell, and the latter's then Rutgers Graduate mate, Shashank Kanade \cite{KR15,KR19,Ru16}.
We are happy to dedicate this note to Doron.

\section{Capparelli's finite sums and their reflections}
\label{sec2}

Before deepening into the topic of the Kanade--Russell identities we first examine a similar but more familiar ground.
We first return to the theme of the Rogers--Ramanujan identity~\eqref{rama1} but now from the perspective of another finite version
\begin{equation*}
\sum_{n\ge0}q^{n^2}\qbin{N-n}{n}
=\sum_{k\in\mathbb Z}(-1)^kq^{k(5k+1)/2}\qbin{N}{\lf(N+5k+1)/2\rf}
\end{equation*}
found by Schur and explicitly stated by Andrews~\cite{An70} (see also \cite[Identity~3.18]{Si03}).
Substituting $1/q$ for $q$ and multiplying the result by $q^{\lf(N/2)^2\rf}$ we obtain, depending on whether $N=2M$ or $N=2M+1$,
\begin{align*}
\sum_{n\ge0}q^{n^2}\qbin{M+n}{2n}
&=\sum_{k\in\mathbb Z}(-1)^kq^{\lf(5k+1)/2\rf^2-k(5k+1)/2}\qbin{2M}{M+\lf(5k+1)/2\rf},
\\
\sum_{n\ge0}q^{n^2+n}\qbin{M+n+1}{2n+1}
&=\sum_{k\in\mathbb Z}(-1)^kq^{\lf5k/2\rf^2+\lf5k/2\rf-k(5k+1)/2}\qbin{2M+1}{M+\lf5k/2\rf+1}
\end{align*}
(we change $n\mapsto M-n$ in the sums on the left-hand sides).
These two (different!) identities are due to Andrews \cite{An81,An86} and listed in \cite[Identities~3.79 and 3.94]{Si03}; their limiting $M\to\infty$ cases
\begin{align*}
\sum_{n\ge0}\frac{q^{n^2}}{(q;q)_{2n}}
&=\frac1{(q;q^2)_\infty(q^4,q^{16};q^{20})_\infty},
\\
\sum_{n\ge0}\frac{q^{n^2+n}}{(q;q)_{2n+1}}
&=\frac1{(q,q^2,q^8,q^9;q^{10})_\infty(q^5,q^6,q^{14},q^{15};q^{20})_\infty}
\end{align*}
appear in Rogers' famous ``second memoir''~\cite{Ro94}.

Secondly, it is natural to draw parallels of the story below with some finite sums from \cite{BU19a,BU19b} underlying Capparelli's identities \cite{Ca95}.
One of this, corresponding to the analytic counterpart \cite[eq.~(7.11)]{KR19}, \cite[Theorem 8]{Ku19b} of the first Capparelli's identity
\begin{equation}
\sum_{m,n\ge0}\frac{q^{2(m^2+3mn+3n^2)}}{(q;q)_m(q^3;q^3)_n}
=(-q^2,-q^4;q^6)_\infty(-q^3;q^3)_\infty,
\label{eq:Cap1}
\end{equation}
is given in~\cite[Theorem~4.3]{BU19b}:
\begin{align}
&
\sum_{m,n\ge0}q^{2(m^2+3mn+3n^2)}\qbin{3N-3m-6n}{m}_q\qbin{2N-2m-3n}{n}_{q^3}
\nonumber\\ &\quad
=\sum_{l\ge0}q^{3(N-2l)(N-2l-1)/2}\qbin{N}{2l}_{q^3}(-q^2,-q^4;q^6)_l.
\nonumber
\end{align}
By considering separately the cases $N=2M$ and $N=2M-1$ and performing the reflection $q\mapsto1/q$ we arrive at
\begin{align*}
&
q^{-6M^2}\sum_{a,b\ge0}q^{2(a^2-3ab+3b^2)}\qbin{3b}{2a}_q\qbin{M+a}{2b}_{q^3}
=q^{-6M^2}\sum_{c\ge0}q^{3c}\qbin{2M}{2c}_{q^3}(-q^2,-q^4;q^6)_{M-c},
\\ &
q^{-6M(M+1)}\sum_{a,b\ge0}q^{2(a^2-3ab+3b^2)+2a-3b-1}\qbin{3b}{2a+1}_q\qbin{M+a}{2b}_{q^3}
\\ &\quad
=q^{-6M(M+1)}\sum_{c\ge0}q^{3c}\qbin{2M}{2c+1}_{q^3}(-q^2,-q^4;q^6)_{M-c},
\end{align*}
respectively, and we find that
\begin{align*}
&
\sum_{a,b\ge0}\frac{q^{2(a^2-3ab+3b^2)}}{(q^3;q^3)_{2b}}\qbin{3b}{2a}
=(-q^2,-q^4;q^6)_\infty\sum_{c\ge0}\frac{q^{3c}}{(q^3;q^3)_{2c}}
\\ &\;
=\frac1{(q^2,q^{10};q^{12})_\infty(q^{3},q^{6},q^{9},q^{9},q^{12},q^{15},q^{15},q^{21},q^{27},q^{33},q^{33},q^{36},q^{39},q^{39},q^{42},q^{45};q^{48})_\infty},
\\ &
\sum_{a,b\ge0}\frac{q^{2(a^2-3ab+3b^2)+2a-3b-1}}{(q^3;q^3)_{2b}}\qbin{3b}{2a+1}
=(-q^2,-q^4;q^6)_\infty\sum_{c\ge0}\frac{q^{3c}}{(q^3;q^3)_{2c+1}}
\\ &\;
=\frac1{(q^2,q^{10};q^{12})_\infty(q^{3},q^{3},q^{9},q^{12},q^{15},q^{18},q^{21},q^{21},q^{27},q^{27},q^{30},q^{33},q^{36},q^{39},q^{45},q^{45};q^{48})_\infty},
\end{align*}
where the right-hand sides were classically summed.

The same identity \eqref{eq:Cap1} admits a different finite version \cite[Theorem~7.1]{BU19a}:
\begin{equation*}
\sum_{m,n\ge0}\frac{q^{2(m^2+3mn+3n^2)}(q^3;q^3)_N}{(q;q)_m(q^3;q^3)_n(q^3;q^3)_{N-m-2n}}
=\sum_{l=-N}^Nq^{l(3l+1)}\qbin{2N}{N+l}_{q^3}.
\end{equation*}
Its $q\mapsto1/q$ reflection after multiplication of both sides by $q^{N(3N+1)}$ reads
\begin{equation*}
\sum_{k,n\ge0}\frac{(-1)^nq^{n(3n+1)/2+k(3N+1)}(q^3;q^3)_N}{(q;q)_{N-k-2n}(q^3;q^3)_k(q^3;q^3)_n}
=\sum_{l=-N}^Nq^{N-l}\qbin{2N}{N+l}_{q^3},
\end{equation*}
and the limit as $N\to\infty$ is uninspiringly equal to $1/(q;q^3)_\infty$.

\section{Finite versions of the Kanade--Russell--Kur\c sung\"oz style double series}
\label{sec3}

The five modulo 9 conjectures about partition generating functions were originally displayed in \cite{KR15,Ru16} through difference equations;
one of these four is already given in the introduction.
Their double-sum Rogers--Ramanujan type versions read
\begin{align}
\label{eq:kr1Sum}
\KR_1(q) &= \sum_{m,n\ge0} \frac{q^{m^2+3mn+3n^2}}{(q;q)_m(q^3;q^3)_n},
\\
\label{eq:kr2Sum}
\KR_2(q) &= \sum_{m,n\ge0} \frac{q^{m^2+3mn+3n^2+m+3n}}{(q;q)_m(q^3;q^3)_n},
\displaybreak[2]\\
\label{eq:kr3Sum}
\KR_3(q) &= \sum_{m,n\ge0} \frac{q^{m^2+3mn+3n^2+2m+3n}}{(q;q)_m(q^3;q^3)_n},
\displaybreak[2]\\
\label{eq:kr4Sum}
\KR_4(q) &= \sum_{m,n\ge0} \frac{q^{m^2+3mn+3n^2+m+2n}}{(q;q)_m(q^3;q^3)_n},
\displaybreak[2]\\
\label{eq:kr5Sum}
\KR_5(q) &= \sum_{m,n\ge0} \frac{q^{m^2+3mn+3n^2+2m+4n}(1+q+q^{m+3n+2})}{(q;q)_m(q^3;q^3)_n}.
\end{align}
Here entries \eqref{eq:kr1Sum}--\eqref{eq:kr4Sum} are found by Kur\c sung\"oz \cite{Ku19} and they correspond to the $I_1$--$I_4$ instances of Kanade--Russell \cite{KR15}, respectively.
The expression \eqref{eq:kr5Sum} corresponds to the later found asymmetric mod~9 conjecture;
its combinatorial version is represented in Russell's thesis \cite{Ru16} and the analytic sum is constructed by us using Kur\c sung\"oz's technique. 
The product sides
\begin{align}
\label{eq:kr1Prod} \KR_1(q) &\overset?= \frac{1}{(q,q^3,q^6,q^8;q^9)_\infty},
\\
\label{eq:kr2Prod} \KR_2(q) &\overset?= \frac{1}{(q^2,q^3,q^6,q^7;q^9)_\infty},
\displaybreak[2]\\
\label{eq:kr3Prod} \KR_3(q) &\overset?= \frac{1}{(q^3,q^4,q^5,q^6;q^9)_\infty},
\displaybreak[2]\\
\label{eq:kr4Prod} \KR_4(q) &\overset?= \frac{1}{(q^2,q^3,q^5,q^8;q^9)_\infty},
\\
\label{eq:kr5Prod} \KR_5(q) &\overset?= \frac{1}{(q,q^4,q^6,q^7;q^9)_\infty}
\end{align}
are precisely the conjectures of Kanade and Russell \cite{KR15,Ru16}.
Notice that \eqref{eq:kr4Prod} is an equivalent form of identity~\eqref{S-lim}. We call the first three products \eqref{eq:kr1Prod}-\eqref{eq:kr3Prod} \textit{symmetric} as for any residue class $i$ mod 9, the residue class $-i$ mod 9 is also present in the product. The remaining products \eqref{eq:kr4Prod} and \eqref{eq:kr5Prod} will then be deemed \textit{asymmetric}.

Following the footsteps of the recent works \cite{BU19a,Un20,Un2} of Uncu, some in collaboration with Berkovich, and using the combinatorial interpretation of the sum sides from \cite{KR15,Ru16} we can write finite generating functions as follows:
\begin{align}
\label{eq:kr1SumN}
\KR_1(q,N)
&= \sum_{m,n\ge0} q^{m^2+3mn+3n^2}
\qbin{N-m-3n+1}{m}_q \qbin{\lfloor\frac23N\rfloor-m-n+1}{n}_{q^3},
\displaybreak[2]\\
\label{eq:kr2SumN}
\KR_2(q,N)
&= \sum_{m,n\ge0} q^{m^2+3mn+3n^2+m+3n}
\qbin{N-m-3n}{m}_q \qbin{\lfloor\frac23N\rfloor-m-n}{n}_{q^3},
\displaybreak[2]\\
\label{eq:kr3SumN}
\KR_3(q,N)
&= \sum_{m,n\ge0} q^{m^2+3mn+3n^2+2m+3n}
\qbin{N-m-3n-1}{m}_q^* \qbin{\lfloor\frac23N\rfloor-m-n}{n}_{q^3},
\displaybreak[2]\\
\label{eq:kr4SumN}
\KR_4(q,N)
&= \sum_{m,n\ge0} q^{m^2+3mn+3n^2+m+2n}
\qbin{N-m-3n}{m}_q \qbin{\lfloor\frac23(N-1)\rfloor-m-n+1}{n}_{q^3},
\displaybreak[2]\\
\KR_5(q,N)
&= \sum_{m,n\ge0} q^{m^2+3mn+3n^2+2m+4n}(1+q)
\nonumber\\ &\qquad\times
\qbin{N-m-3n-1}{m}_q \qbin{\lfloor\frac23(N-1)\rfloor-m-n+\delta_{3\mid(N-2)}}{n}_{q^3}
\nonumber\\
\label{eq:kr5SumN}
&\kern-20mm
+ \sum_{m,n\ge0} q^{m^2+3mn+3n^2+3m+7n+2}
\qbin{N-m-3n-2}{m}_q \qbin{\lfloor\frac23(N-1)\rfloor-m-n}{n}_{q^3},
\end{align}
where the asterisk in \eqref{eq:kr3SumN} means that the $q$-binomial $\qbin{N-m-3n-1}{m}_q$ is understood as~1 when it becomes $\qbin{-1}{0}_q$ (in other words, when $m=N-m-3n=0$, which may only occur when $N\equiv0\pmod3$);
also $\delta_{a\mid b}$ used in~\eqref{eq:kr5SumN} stands for 1 if $a\mid b$ and for 0 otherwise.
Then clearly the limiting $N\to\infty$ cases of \eqref{eq:kr1SumN}--\eqref{eq:kr5SumN} become \eqref{eq:kr1Sum}--\eqref{eq:kr5Sum};
it is also routine to verify that the polynomial (finite) versions of the latter double sums indeed satisfy the recurrence equations given in \cite{KR15,Ru16}.
Since we find the technique of converting infinite sums into finite versions important, we select it for preservation in Appendix~\ref{appA} below.

In all the expressions \eqref{eq:kr1Sum}--\eqref{eq:kr5Sum} the exponent of $q$ keeps track of the size of the counted partitions.
In fact, the technique developed in \cite{BU19a,Un20,Un2, Ku19} allows us to write the generating functions with an additional statistics by multiplying the $(m,n)$th term in the sum by $x^{2n+m}$ in \eqref{eq:kr1SumN}--\eqref{eq:kr4SumN} and taking, for instance
\begin{align*}
\KR_5(q,N;x)
&= \sum_{m,n\ge0} q^{m^2+3mn+3n^2+2m+4n}(1+xq)
\qbin{N-m-3n-1}{m}_q
\\ &\qquad\times
\qbin{\lfloor\frac23(N-1)\rfloor-m-n+\delta_{3\mid(N-2)}}{n}_{q^3} x^{m+2n}
\displaybreak[2]\\
&\quad
+ \sum_{m,n\ge0} q^{m^2+3mn+3n^2+3m+7n+2}
\qbin{N-m-3n-2}{m}_q
\\ &\qquad\times
\qbin{\lfloor\frac23(N-1)\rfloor-m-n}{n}_{q^3} x^{m+2n+1},
\end{align*}
Then the exponent of $x$ in these general generated functions keeps track of the number of parts in the counted partitions.
We leave this as a side remark but keep a hope that the more involved series can be used for approaching conjectures~\eqref{eq:kr1Prod}--\eqref{eq:kr5Prod}.

\section{Modulo 45 reflections}
\label{sec4}

In our language here, Warnaar's conjectures from the introduction can be stated as follows.
{\em Define
\begin{align*}
\RK_4(q,3M)&=q^{M(3M+2)} \KR_4(1/q,3M), \\
\RK_4(q,3M+1)&=q^{M(3M+5)} \KR_4(1/q,3M+1), \\
\RK_4(q,3M+2)&=q^{(M+1)(3M+2)} \KR_4(1/q,3M+2).
\end{align*}
Then}
\begin{align*}
\RK_4(q,3\infty)
=\lim_{M\to\infty} \RK_4(q,3M)
&\overset?= \frac{1}{(q^2;q^3)_\infty (q^3,q^9,q^{12},q^{21},q^{30},q^{36},q^{39};q^{45})_\infty},
\\
\RK_4(q,3\infty+1)
=\lim_{M\to\infty} \RK_4(q,3M+1)
&\overset?= \frac{1}{(q^2;q^3)_\infty (q^3,q^{12},q^{18},q^{21},q^{27},q^{30},q^{39};q^{45})_\infty},
\\
\RK_4(q,3\infty+2)
=\lim_{M\to\infty} \RK_4(q,3M+2)
&\overset?=\RK_4(q,3\infty)+q^2\RK_4(q,3\infty+1).
\end{align*}

To convert the limits into infinite sums we use the finite sum representation \eqref{eq:kr4SumN}.
We deduce that
\begin{align*}
\RK_4(q,3M)
&=q^{M(3M+2)} \KR_4(1/q,3M)
\\
&=q^{3M^2+2M}\sum_{m,n\ge 0} q^{m^2+3mn +3n^2-(2n+m)(3M+1)} \qbin{3M-3n-m}{m}_q \qbin{2M-n-m}{n}_{q^3}.
\end{align*}
Changing the summation to the one over $a=3M-2m-3n$, $b=2M-m-2n$ (equivalently, $m=3b-2a$, $n=M+a-2b$) the sum transforms into
\[
\RK_4(q,3M)
=\sum_{a,b\ge 0} q^{a^2 - 3ab + 3b^2 + b} \qbin{3b - a}{a}_q \qbin{M+a-b}{b}_{q^3}
\]
leading to
\begin{align}
\RK_4(q,3\infty) 
&= \sum_{a,b\ge 0} \frac{q^{a^2 - 3ab + 3b^2 + b}}{(q^3;q^3)_b} \qbin{3b - a}{a}
\nonumber\\
&\overset?= \frac{1}{(q^2;q^3)_\infty (q^3,q^9,q^{12},q^{21},q^{30},q^{36},q^{39};q^{45})_\infty} .
\label{eq:KR4_reflected_0mod3}
\end{align}
Similarly, we have 
\begin{align}
\RK_4(q,3M+1)
&= \sum_{a,b\ge 0} q^{a^2 - 3ab + 3b^2 + b-2} \qbin{3b - a - 1}{a}_q \qbin{M+a-b+1}{b}_{q^3},
\nonumber\\
\RK_4(q,3M+2)
&= \sum_{a,b\ge 0} q^{a^2 - 3ab + 3b^2 + b} \qbin{3b - a + 1}{a}_q \qbin{M+a-b}{b}_{q^3},
\nonumber\\ \intertext{hence}
\RK_4(q,3\infty+1)
&= \sum_{a,b\ge 0} \frac{q^{a^2 - 3ab + 3b^2 + b-2}}{(q^3;q^3)_b} \qbin{3b - a - 1}{a}
\nonumber\\
&\overset?= \frac{1}{(q^2;q^3)_\infty (q^3,q^{12},q^{18},q^{21},q^{27},q^{30},q^{39};q^{45})_\infty},
\label{eq:KR4_reflected_1mod3}
\displaybreak[2]\\
\RK_4(q,3\infty+2)
&= \sum_{a,b\ge 0} \frac{q^{a^2 - 3ab + 3b^2 + b}}{(q^3;q^3)_b} \qbin{3b - a+1}{a}
\nonumber\\
&=\RK_4(q,3\infty)+q^2\RK_4(q,3\infty+1).
\label{eq:KR4_reflected_2mod3}
\end{align}
Here the last part \eqref{eq:KR4_reflected_2mod3} follows from
\begin{equation*}
\sum_{a,b\ge 0} \frac{q^{a^2 - 3ab + 3b^2 + b}}{(q^3;q^3)_b}
\bigg(\qbin{3b - a - 1}{a}+\qbin{3b - a}{a}-\qbin{3b - a + 1}{a}\bigg)
=0,
\end{equation*}
which is a consequence of simple single-sum evaluation
\begin{equation}
\sum_{a=0}^{L+1} q^{a(a-L)}
\bigg(\qbin{L - a - 1}{a}+\qbin{L - a}{a}-\qbin{L - a + 1}{a}\bigg)
=0
\quad\text{for}\; L=0,1,2,\dotsc.
\label{eq:id}
\end{equation}

Quite remarkably a similar structure is inherited by the other `asymmetric' conjecture~\eqref{eq:kr5Prod}.
{\em Define
\begin{align*}
\RK_5(q,3M)&=q^{M(3M+1)} \KR_5(1/q,3M), \\
\RK_5(q,3M+1)&=q^{(M+1)(3M+1)} \KR_5(1/q,3M+1), \\
\RK_5(q,3M+2)&=q^{(M+2)(3M+1)} \KR_5(1/q,3M+2).
\end{align*}
Then}
\begin{align*}
\lim_{M\to\infty} \RK_5(q,3M+1)
&\overset?= \frac{1}{(q;q^3)_\infty (q^6,q^{9},q^{15},q^{24},q^{33},q^{36},q^{42};q^{45})_\infty},
\\
\lim_{M\to\infty} \RK_5(q,3M+2)
&\overset?= \frac{1}{(q;q^3)_\infty (q^6,q^{15},q^{18},q^{24},q^{27},q^{33},q^{42};q^{45})_\infty},
\\
\lim_{M\to\infty}\RK_5(q,3M)
&\overset?=\lim_{M\to\infty} \RK_5(q,3M+1)
+q^2 \lim_{M\to\infty} \RK_5(q,3M+2).
\end{align*}
Our analysis as before translates the expectations into
\begin{align}
\RK_5(q,3\infty+1)
&= \sum_{a,b\ge 0} \frac{ q^{a^2- 3 a b + 3 b^2 + 2 b } (1 + q )}{(q^3;q^3)_b} \qbin{3 b - a}{a}
\nonumber\\ &\quad
+\sum_{a,b\ge 0} \frac{q^{a^2- 3 a b + 3 b^2 -a + 5 b+2}}{(q^3;q^3)_b} \qbin{3 b - a + 1}{a}
\nonumber\\
&\overset?= \frac{1}{(q;q^3)_\infty (q^6,q^{9},q^{15},q^{24},q^{33},q^{36},q^{42};q^{45})_\infty},
\label{eq:KR5_reflected_1mod3}
\displaybreak[2]\\
\RK_5(q,3\infty+2)
&=\sum_{a,b\ge 0}  \frac{q^{a^2- 3 a b + 3 b^2 + 2 b-2 } (1 + q )}{(q^3;q^3)_b}\qbin{3 b - a - 1}{a}
\nonumber\\ &\quad
+\sum_{a,b\ge 0} \frac{q^{a^2- 3 a b + 3 b^2 -a+ 5 b }}{(q^3;q^3)_b} \qbin{3 b - a}{a}
\nonumber\\
&\overset?= \frac{1}{(q;q^3)_\infty (q^6,q^{15},q^{18},q^{24},q^{27},q^{33},q^{42};q^{45})_\infty},
\label{eq:KR5_reflected_2mod3}
\displaybreak[2]\\
\RK_5(q,3\infty)
&= \sum_{a,b\ge 0} \frac{ q^{a^2- 3 a b + 3 b^2 + 2 b } (1 + q )}{(q^3;q^3)_b} \qbin{3 b - a+1}{a}
\nonumber\\ &\quad
+\sum_{a,b\ge 0} \frac{q^{a^2- 3 a b + 3 b^2 -a + 5 b+2}}{(q^3;q^3)_b} \qbin{3 b - a + 2}{a}
\nonumber\\
&=\RK_5(q,3\infty+1)+q^2\RK_5(q,3\infty+2).
\label{eq:KR5_reflected_0mod3}
\end{align}
Again, the equality in \eqref{eq:KR5_reflected_0mod3} follows from \eqref{eq:id}.

In spite of the suggested simplicity of the reflections in the asymmetric cases, our investigation of the symmetric ones \eqref{eq:kr1Prod}--\eqref{eq:kr3Prod} brings to life somewhat more sophisticated expectations:
\begin{align}
&
\sum_{a,b\ge 0} \frac{q^{a^2-3ab+3b^2-1}}{(q^3;q^3)_b}\qbin{3b-a-1}{a}_q
=\lim_{M\to\infty}q^{3M(M+1)}\KR_1(q^{-1},3M)
=\RK_1(q,3\infty)
\nonumber\\ &\;
\overset?=\<2,8,11,20\>+q^3\<2,14,20,22\>-q^8\<17,19,20,22\>
\nonumber\\ &\;
=\<1,8,13,20\>-q\<4,7,13,20\>+q^5\<7,16,17,20\>,
\label{eq:KR1_reflected_0mod3}
\displaybreak[2]\\
&
\sum_{a,b\ge 0} \frac{q^{a^2-3ab+3b^2}}{(q^3;q^3)_b}\qbin{3b-a}{a}_q
=\lim_{M\to\infty}q^{3(M+1)^2}\KR_1(q^{-1},3M+2)
=\RK_1(q,3\infty+2)
\nonumber\\ &\;
\overset?=\<1,7,11,20\>+q^6\<11,13,14,20\>-q^6\<8,14,19,20\>
\nonumber\\ &\;
=\<1,4,17,20\>-q^4\<2,16,19,20\>-q^5\<4,16,20,22\>,
\label{eq:KR1_reflected_2mod3}
\displaybreak[2]\\
&
\sum_{a,b\ge 0} \frac{q^{a^2-3ab+3b^2-a+3b}}{(q^3;q^3)_b}\qbin{3b-a}{a}_q
=\lim_{M\to\infty}q^{3M(M+1)}\KR_2(q^{-1},3M)
=\RK_2(q,3\infty)
\nonumber\\ &\;
\overset?=\<2,5,14,22\>-q^2\<5,7,16,17\>-q^5\<5,17,19,22\>
\nonumber\\ &\;
=q^{-3}\<1,5,8,13\>-q^{-3}\<2,5,8,11\>-q^{-2}\<4,5,7,13\>,
\label{eq:KR2_reflected_0mod3}
\displaybreak[2]\\
&
\sum_{a,b\ge 0} \frac{q^{a^2-3ab+3b^2-a+3b}}{(q^3;q^3)_b}\qbin{3b-a+1}{a}_q
=\lim_{M\to\infty}q^{3(M+1)^2-1}\KR_2(q^{-1},3M+2)
=\RK_2(q,3\infty+2)
\nonumber\\ &\;
\overset?=\<2,5,16,19\>-q^2\<5,8,14,19\>+q^2\<5,11,13,14\>
\nonumber\\ &\;
=q^{-4}\<1,4,5,17\>-q^{-4}\<1,5,7,11\>-q\<4,5,16,22\>,
\label{eq:KR2_reflected_2mod3}
\displaybreak[2]\\
&
1+\sum_{a,b\ge 0} \frac{q^{a^2-3ab+3b^2+a}}{(q^3;q^3)_b}\qbin{3b-a-2}{a}_q
=\lim_{M\to\infty}q^{3M(M+1)}\KR_3(q^{-1},3M)
=\RK_3(q,3\infty)
\nonumber\\ &\;
\overset?=\<4,7,10,13\>-q^4\<7,10,16,17\>-q^7\<10,17,19,22\>
\nonumber\\ &\;
=q^{-1}\<1,8,10,13\>-q^{-1}\<2,8,10,11\>-q^2\<2,10,14,22\>,
\label{eq:KR3_reflected_0mod3}
\displaybreak[2]\\
&
\sum_{a,b\ge 0} \frac{q^{a^2-3ab+3b^2+a-2}}{(q^3;q^3)_b}\qbin{3b-a-1}{a}_q
=\lim_{M\to\infty}q^{3M(M+2)}\KR_3(q^{-1},3M+2)
=\RK_3(q,3\infty+2)
\nonumber\\ &\;
\overset?=\<2,10,16,19\>+q\<4,10,16,22\>-q^2\<8,10,14,19\>
\nonumber\\ &\;
=q^{-4}\<1,4,10,17\>-q^{-4}\<1,7,10,11\>-q^2\<10,11,13,14\>,
\label{eq:KR3_reflected_2mod3}
\end{align}
where the summands are (modular) products
\begin{equation}
\<c_1,c_2,c_3,c_4\>=\frac{(q^{45};q^{45})_\infty}{(q^3;q^3)_\infty\prod_{j=1}^4(q^{c_j},q^{45-c_j};q^{45})_\infty}.
\label{eq:c}
\end{equation}
Furthermore, identity \eqref{eq:id} and its modifications also give us (provably!)
\begin{align}
\sum_{a,b\ge 0} \frac{q^{a^2-3ab+3b^2}}{(q^3;q^3)_b}\qbin{3b-a+1}{a}_q
&=\lim_{M\to\infty}q^{3M(M+1)+1}\KR_1(q^{-1},3M+1)
\nonumber\\ &\kern-24mm
=\RK_1(q,3\infty+1)
=q\RK_1(q,3\infty)+\RK_1(q,3\infty+2),
\label{eq:KR1_reflected_1mod3}
\displaybreak[2]\\
\sum_{a,b\ge 0} \frac{q^{a^2-3ab+3b^2-a+3b}}{(q^3;q^3)_b}\qbin{3b-a+2}{a}_q
&=\lim_{M\to\infty}q^{3M(M+1)}\KR_2(q^{-1},3M+1)
\nonumber\\ &\kern-24mm
=\RK_2(q,3\infty+1)
=\RK_2(q,3\infty)+\RK_2(q,3\infty+2),
\label{eq:KR2_reflected_1mod3}
\displaybreak[2]\\
\sum_{a,b\ge 0} \frac{q^{a^2-3ab+3b^2+a}}{(q^3;q^3)_b}\qbin{3b-a}{a}_q
&=\lim_{M\to\infty}q^{3M(M+1)}\KR_3(q^{-1},3M+1)
\nonumber\\ &\kern-24mm
=\RK_3(q,3\infty+1)
=\RK_2(q,3\infty)+q^2\RK_2(q,3\infty+2).
\label{eq:KR3_reflected_1mod3}
\end{align}
One interesting outcome of these conjectures is the positivity of the $q$-expansions of the linear combinations of products that appear on the right-hand sides for all $\RK_i(q,3\infty+l)$ for $i=1,2,3$ and $l=0,1,2$.
In fact, we have checked numerically that the combinations remain positive after each of them is multiplied by $(q^3;q^3)_\infty/(q^{45};q^{45})_\infty$, but for that we have no explanation.

The complexity of the reflections serves a good reason for lack of simple single-sum evaluations for the original sums $\KR_i(q,N)$, where $i=1,2,3$.

One may also argue that the (combinatorially motivated!) finite versions \eqref{eq:kr1SumN}--\eqref{eq:kr5SumN} may be not best in approaching the Kanade--Russell identities \eqref{eq:kr1Prod}--\eqref{eq:kr5Prod} and one can possibly try more general sums like
\[
F(N,M):=\sum_{m,n\ge0} q^{m^2+3mn+3n^2}
\qbin{N-m-3n}{m}_q \qbin{M-m-n}{n}_{q^3}
\]
for different specializations $M=M(N)\to\infty$ as $N\to\infty$ instead. Although using \cite{qFunctions} it is easy to observe that $F(N,M)$ satisfy simple-looking recurrences such as
\[
F(N,M) = F(N-1,M) +  q^{N-1} F(N-2, M-1) \quad\text{for}\; 3\nmid N ,
\]
none of these variations seem to produce (linear combinations of) products for the limits of the sums reflected under $q\mapsto1/q$.
In other words, the combinatorics of the Kanade--Russell partition functions is a delicate sensor of arithmetic features.

\section{Modular remarks}
\label{sec5}

Without giving precise definitions of modular (and mock modular) functions,
for which the reader can consult some standard sources (e.g., \cite{On09,Za09}),
we notice that it was the modularity of the product sides in \eqref{eq:kr1Prod}--\eqref{eq:kr3Prod} that led us to reasonable guesses for \eqref{eq:KR1_reflected_0mod3}--\eqref{eq:KR3_reflected_1mod3}.
The use of the servers at the Research Institute for Symbolic Computation (RISC) for the actual calculations significantly cooled down our laptops.

The products in \eqref{eq:kr4Prod}, \eqref{eq:kr5Prod} are not modular but, expectedly, mock modular.
In fact quite likely, they form a 2-vector that exhibits a modular behavior,
and similar may be formed from the products in \eqref{eq:KR4_reflected_0mod3}, \eqref{eq:KR5_reflected_1mod3}
and in \eqref{eq:KR4_reflected_1mod3}, \eqref{eq:KR5_reflected_2mod3}.

The main obstacle for using the powerful (mock) modular structure in proving the Kanade--Russell conjectures and their reflected counterparts is a difficulty in establishing the modular behavior for the sum sides of the identities.

\medskip
\noindent
\textbf{Acknowledgements.}
We are greatly intended to Ole Warnaar for sharing with us his $\KR_4$ observations and for his patience in trivializing many intermediate observations by these authors, like explaining that KR-looking
\begin{align*}
\sum_{m,n\ge0}\frac{q^{m^2+2mn+2n^2}}{(q;q)_m(q^2;q^2)_n}
&=\frac{(q^3;q^3)_\infty^2}{(q;q)_\infty(q^6;q^6)_\infty},
\\
\sum_{m,n\ge0}\frac{q^{m^2+2mn+2n^2+m+2n}}{(q;q)_m(q^2;q^2)_n}
&=\frac{(q^6;q^6)_\infty^2}{(q^2;q^2)_\infty(q^3;q^3)_\infty}
\end{align*} 
are instances of Bressoud's identities from~\cite{Br79}.

The initial discussions about this project took place during the second author's visit in the Research Institute for Symbolic Computation (Linz) in February 2020.
We thank Peter Paule for making that visit not only possible but also mathematically fruitful.
The conversations with colleagues at RISC and RICAM in Linz were another source of inspiration.
After the trip this project went entirely online (together with the academic life more generally) with a small live break for the \emph{85th S\'eminaire Lotharingien de Combinatoire} in Strobl in September 2020.
We thank Christian Krattenthaler for organizing the event and exchange of mathematical ideas.
Finally, we are happy to thank Stepan Konenkov and Eric Mortenson for catching typos in the final version of this text.

The first author would like to extend gratitude to UK Research and Innovation EPSRC and Austrian Science Fund (FWF) for supporting his research through the grants EP/T015713/1 and FWF P-34501-N, respectively.

\appendix
\section{Combinatorial construction of the finite versions of the Kanade--Russell--Kur\c sung\"oz style double series}
\label{appA}

We give a brief description and considerations that go into identifying the polynomial analogues \eqref{eq:kr1SumN}--\eqref{eq:kr5SumN} of \eqref{eq:kr1Sum}-\eqref{eq:kr5Sum}.
To that end, we present how we \emph{combinatorially} construct the refinement of \eqref{eq:kr4Sum}.
Recall that \eqref{eq:kr4Sum} was shown by Kur\c{s}ung\"oz \cite{Ku19} to be the generating function for the number of partitions that satisfy the gap conditions prescribed in the $I_4$ conjecture of \cite{KR15}. Here we prove that the finite analogue \eqref{eq:kr4SumN} is the generating function for the number of partitions into parts $\le N$ that satisfy the gap conditions prescribed in the $I_4$ conjecture\,---\,the conditions are given as \eqref{cond-a}--\eqref{cond-d} in the Introduction.
We will be mimicking the constructions of \cite{BU19a,Un20,Un2}, and invite the interested reader to examine these references to see longer expositions of this technique. 

Recall conditions \eqref{cond-b}, \eqref{cond-c} from the second page.
The $I_4$ conjecture considers partitions in which 1 does not appear as a part, the difference between parts is $\ge3$ at distance 2 such that if two successive parts differ by $\le1$, then their sum is congruent to $2\pmod3$.
First we would like to explain how to interpret, partition-theoretically, the pieces of \eqref{eq:kr4Sum},
\[
\sum_{m,n\ge0 } \frac{q^{m^2+3mn+3n^2+m+2n}}{(q;q)_m(q^3;q^3)_n},
\]
as the generating function for these partitions.
This is done in the spirit of Kur\c sung\"oz's construction in \cite{Ku19} using the vocabulary of \cite{BU19a,Un20,Un2}.
For the rest of this discussion, we consider a partition to be a finite sequence of \textit{non-decreasing} positive integers.
The $q$-factor $q^{m^2+3mn+3n^2+m+2n}$ is the size of the partition
\[
\pi_{m,n}=(\underline{2,3},\underline{5,6},\dots,\underline{3n-1,3n}, 3n+2, 3n+4,\dots,3n+2m).
\]
We call the underlined consecutive parts of $\pi_{m,n}$ \textit{pairs} and the rest of the terms \textit{singletons}.
Notice that $\pi_{m,n}$ is a partition that satisfies the gap conditions of~$I_4$.
Moreover, it is not hard to observe that $\pi_{m,n}$ is the partition with the smallest possible size that satisfies the gap conditions of $I_4$ into $2n+m$ parts,
where the minimal distance condition
(if two successive parts differ by $\le1$, then their sum is congruent to $2\pmod3$)
appears exactly $n$~times.
We call $\pi_{m,n}$ a \textit{minimal configuration} with $2n+m$ parts and $n$ minimal gaps.

If one adds any non-negative integer value $r_m$ to the largest part $3n+2m$ of $\pi_{m,n}$,
the outcome partition still satisfies the gap conditions of~$I_4$.
Similarly, after adding $r_m$ to the largest part, if one adds some non-negative integer $r_{m-1}\le r_m$ to the second largest value of $\pi_{m,n}$, the outcome partition satisfies the gap conditions of~$I_4$.
Repeating this process, we can conclude that if one adds $(r_1,\dots, r_m)$, where $0\le r_1\le r_2\le\dots\le r_m$, to the $m$ largest parts (the singletons) of $\pi_{m,n}$, the outcome partition $\pi_{n}$ still satisfies the gap conditions of~$I_4$.
Moreover, this addition can be reversed.
Given a partition $\pi_{n}$, where the singletons are possibly not at their original locations while the pairs are left untouched, we can easily recover the list $(r_1,\dots, r_m)$ and the minimal configuration $\pi_{n,m}$ that gives rise to~$\pi_{n}$.
Finally, we note that the lists $(r_1,\dots, r_m)$, where $0\le r_1\le\dots\le r_m$, are in one-to-one correspondence to partitions into $\le m$ parts.
The generating function for the partitions into $\le m$ parts is $(q;q)^{-1}_m$.
Hence, now we can conclude that
\[
\frac{q^{m^2+3mn+3n^2+m+2n}}{(q;q)_m}
\]
is the generating function for the number of all the partitions of the form
\[
\pi_{n}=(\underline{2,3},\underline{5,6},\dots,\underline{3n-1,3n}, s_1, s_2,\dots, s_m),
\]
where $3n+2\le s_1$ and $s_{i}-s_{i-1}\ge2$ for $i=2,\dots,m$.

Now, given a partition $\pi_n$ with $m$ singletons, assuming that there are no close-by singletons, we can move the largest pair to the next permissible pair and repeat this (reversible) process to generate other partitions with $n$~arbitrary pairs and $m$~singletons.
To be precise, there are two possible (free) forward motions of pairs and these are
\[
\underline{3k-1,3k}\mapsto\underline{3k+1,3k+1}
\quad\text{and}\quad
\underline{3k+1,3k+1}\mapsto\underline{3k+2,3k+3}.
\]
In both, forward motions of a pair the total size change of the partition is~3.
Similar to moving singletons, by moving the largest pair, followed by the second largest pair, \emph{etc}., we never need to worry about pairs crossing each other.
However, although the forwards motions of singletons were done freely, as we move a pair, we may come close to a singleton and violate the gap conditions of~$I_4$.
To avoid this, we define (reversible) crossing-over rules of pairs over singletons.
The two necessary crossing-over rules are as follows: 
\begin{align*}
\underline{3k-1,3k},3k+2&\mapsto 3k-1,\underline{3k+2,3k+3}, \\
\underline{3k+1,3k+1},3k+4&\mapsto 3k+1,\underline{3k+4,3k+4}.
\end{align*}
Notice that these motions also add 3 to the total size of the partition.
Hence, now given a partition $\pi_n$, starting from the largest pair, we can move the pairs forward and make new partitions into $2n+m$ parts that has $n$~pairs which satisfy the gap conditions of~$I_4$.
Similar to the singletons' case, the forwards motion lists of the pairs corresponds to partitions into $\le n$ parts.
Since each motion of the pairs add 3 to the total size of the partition,
we instead use $(q^3;q^3)^{-1}_\infty$ here.
Therefore, the summand of \eqref{eq:kr4Sum},
\[
\frac{q^{m^2+3mn+3n^2+m+2n}}{(q;q)_m(q^3;q^3)_n},
\]
is the generating function of all the partitions that satisfy the gap conditions of $I_4$ which has $n$~pairs and $m$~singletons.
Summing over all the $m$ and $n$ gives the generating function of all the partitions that satisfy the gap conditions of~$I_4$ and finishes the combinatorial construction/interpretation of~\eqref{eq:kr4Sum} in the spirit of Kur\c sung\"oz. 

For the finite analogue, all we need to do is to restrict the forward motions of the terms.
If we want to add the new restriction that all the parts of our partitions are $\le N$ (condition~\eqref{cond-a} on the second page), the forwards motion of the largest singleton cannot be free.
In this situation $r_m$ must be $\le N-(3n+2m)$, hence the generating function for partitions into $\le m$ parts, which represented the reversible forward motion of the singletons, is now replaced by the generating function for the partitions into $\le m$ parts with each part $\le N-(3n+2m)$.
This generating function is
\[
\qbin{N-(3n+2m)+m}{m}_q.
\]
Similarly, we need to adjust the reversible forwards motion of the pairs.
At each move of a pair, the middle point of a pair (the arithmetic mean of the elements in the pair) moves $2/3$ steps.
From the largest pair $\underline{3n-1,3n}$ to the upper bound $N$, there are exactly $\lfloor 2(N-3n )/3 \rfloor + \delta_{3\mid (N-1)}$ many steps that a pair can move forwards if all the motions are done via free motion rules.
As a pair moves forward, it can cross up to $m$~singletons, and crossing over each singleton also means missing one possible location where the pair could have stopped.
Hence, the actual number of steps a pair can move forward is $\lfloor 2(N-3n )/3 \rfloor + \delta_{3\mid (N-1)} - m$.
Therefore, we need to replace the generating function $(q^3;q^3)^{-1}_n$ with the $q$-binomial coefficient
\[
\qbin{\lfloor 2(N-3n )/3 \rfloor + \delta_{3\mid (N-1)} - m + n}{n}_{q^3}
\]
to restrict the forward motion of the $n$~pairs, where all the parts still remain $\le N$ after the motions.
Compiling the above discussion we get that
\[
\sum_{m,n\ge0} q^{m^2+3mn+3n^2+m+2n}
\qbin{N-m-3n}{m}_q \qbin{\lfloor 2(N-3n )/3 \rfloor + \delta_{3\mid (N-1)} - m + n }{n}_{q^3}
\]
is the generating function for the partitions into parts $\le N $ which satisfy the gap conditions of~$I_4$,
and this expression is equal to \eqref{eq:kr4SumN} after simplification of the top argument of the second $q$-binomial coefficient.


\end{document}